\documentclass{amsart}
\usepackage{amssymb}
\newtheorem{theorem}{Theorem}[section]
\newtheorem{lemma}[theorem]{Lemma}
\newtheorem{proposition}[theorem]{Proposition}

\usepackage[english]{babel} 
\theoremstyle{definition}

\theoremstyle{remark}
\newtheorem{remark}[theorem]{Remark}

\numberwithin{equation}{section}

\renewcommand{\dim}{\mathrm{dim}}

\renewcommand{\span}{\mathrm{span}}

\newcommand{\conv}{\mathrm{conv}}

\newcommand{\R}{\mathbb{R}}
\newcommand{\N}{\mathbb{N}}

\newcommand{\X}{\mathrm{X}}
\renewcommand{\H}{\mathrm{H}}
\newcommand{\Y}{\mathrm{Y}}
\newcommand{\Z}{\mathrm{Z}}
\newcommand{\B}{\mathbf{B}}
\newcommand{\I}{\mathbf{I}}

\renewcommand{\S}{\mathbf{S}}

\renewcommand{\ker}{\mathrm{Ker}}

\begin{document}

\title{Directionally Euclidean Structures of Banach Spaces}

\begin{abstract}
We study spaces with directionally asymptotically controlled ellipsoids approximating the unit ball in finite-dimensions. These ellipsoids are the unique minimum volume ellipsoids, which contain the unit ball of the corresponding 
finite-dimensional subspace. The directional control here means that we evaluate the ellipsoids with a given functional of the dual space. The term asymptotical refers to the fact that we take '$\limsup$' over finite-dimensional subspaces. 

This leads to some isomorphic and isometric characterizations of Hilbert spaces. An application involving Mazur's rotation problem is given. We also discuss the complexity of the family of ellipsoids as the dimension and geometry varies. 
\end{abstract}

\author{Jarno Talponen}
\address{Aalto University, Institute of Mathematics, P.O. Box 11100, FI-00076 Aalto, Finland} 
\email{talponen@cc.hut.fi}
\date{\today}
\maketitle

\section{Introduction}

This paper deals with the local theory of Banach spaces. Here that roughly means that we will make deductions about the geometry of Banach spaces by 
studying the asymptotical behaviour of its finite-dimensional subspaces as the dimension grows. 

In a sense, there is a natural way of approximating a norm of a finite-dimensional normed space by a norm induced by an inner product. 
Namely, by compactness there exists an ellipsoid with the minimal volume (i.e. minimal Lebesgue measure) such that the ellipsoid contains the unit ball.
It was an interesting observation made by Auerbach that such a minimum volume ellipsoid is unique and thus preserved in linear isometries of the normed space (see \cite{Auerbach}, \cite{AuMU}). 
Thus, by now it is a rather standard practice to use the Hilbertian norm induced by the minimum volume ellipsoid to approximate the original norm. Although there are several other approaches as well, see \cite{GM}.

In Hilbert spaces the minimum volume ellipsoids in finite-dimensional subspaces are precisely the unit balls of the corresponding subspaces and a fortiori uniformly bounded. On the other and, 
it is not difficult to see that if the Banach space in question is not isomorphic to a Hilbert space, then the diameter of the minimum volume ellipsoids is not uniformly bounded. 
It turns out that the diameter of the ellipsoids may not be uniformly bounded even if the space is isomorphically Hilbertian. 
Here we will study spaces such that the ellipsoids are asymptotically bounded in some direction. More precisely, the term 'Directionally Euclidean Structures' in the title refers to spaces $\X$ with $f\in \X^{\ast},\ f\neq 0,$ such that 
that for sufficiently large finite-dimensional subspaces $F$ the corresponding minimum volume ellipsoid $\mathcal{E}_{F}$ has the property that its image under $|f|$ is bounded by a uniform constant. 
This will be stated accurately shortly.

To make a short introductory comment about the content of this paper and its connections, this topic has a flavour slightly similar to the study of weak Hilbert spaces and it relies heavily on ultrafilter based analysis.
The methodology and the problem setting here are closely related to the papers \cite{Ca} and \cite{T}. 

We will characterize Banach spaces isomorphic to Hilbert spaces as those spaces, which have directionally Euclidean structure in every direction, up to isomorphism. 
As an application, we will obtain a result related to Mazur's rotation problem. 

\subsection{Preliminaries}
We refer to \cite{HHZ}, \cite{Heinrich}  and \cite{S} for suitable general background information. The known 'Mazur's rotation problem' appearing in Banach's book is discussed extensively in the survey \cite{BR} 
and for local theory of Banach spaces involving ellipsoids we refer to \cite{Pisier}. 

Here $\X$ and $\Y$ stand for real Banach spaces, $\B_{\X}$ is the closed unit ball and $\S_{\X}$ the unit sphere of $\X$. 
We denote by $\mathrm{Aut}(\X)$ the group of isomorphisms $T\colon \X\rightarrow \X$. The rotation group $\mathcal{G}_{\X}$ is the subgroup of all the isometries in $\mathrm{Aut}(\X)$. The identity map $I\colon \X\to \X$ is the neutral element. 
The group of finite-dimensional perturbations of the identity is given by 
\[\mathcal{G}_{F}=\{T\in \mathcal{G}_{\X}:\ \mathrm{Rank}(I-T)<\infty\}.\]
We say that $\X$ is \emph{convex-transitive} with respect to $\mathcal{G}\subset\mathcal{G}_{\X}$ if for all $x\in \S_{\X}$ it holds that 
$\overline{\conv}(\{Tx:\ T\in \mathcal{G}\})=\B_{\X}$. A stronger condition than convex-transitivity is \emph{almost transitivity} w.r.t. $\mathcal{G}$, namely, that $\overline{\{Tx:\ T\in \mathcal{G}\}}=\S_{\X}$ for $x\in\S_{\X}$.

We call a bilinear map $B\colon \X\times \X\to \R$ such that $B(x,x)\geq 0$ for $x\in \X$ and $B(z,z)=0$ for some $z\neq 0$, a \emph{degenerate inner product}. If additionally $B(x,x)>0$ for all $x\neq 0$, then $B$ is an inner product. 
Sometimes we also call $B$, slightly overemphasising, a non-degenerate inner product. An ellipsoid is a set of the form $\{x\in E:\ (x|x)\leq 1\}$, where $(\cdot|\cdot)$ is an inner product on a finite-dimensional (sub)space $E$. 
The set of finite-dimensional subspaces of $\X$ will usually be denoted by $\mathcal{F}$.

\section{Directionally bounded ellipsoids}

For a finite-dimensional subspace $E\subset \X$ and $f\in \S_{\X^{\ast}}$ put
\[\eta(f,E)=\max\{|fy|:\ y\in \mathcal{E}_{E}\},\]
where $\mathcal{E}_{E}\subset E$ is the unique minimum volume ellipsoid containing $\B_{\X}\cap E$.
We say that a Banach space $\X$ has \emph{Directionally Euclidean Structure} (DES) in direction $f\in \S_{\X^{\ast}}$ if 
\begin{equation}\label{eq: LES}
\inf_{F}\sup_{E}\eta(f,E)<\infty,
\end{equation}
where the infimum is taken over finite-dimensional spaces $F\subset\X$ and the supremum is taken over finite-dimensional subspaces 
$E\subset \X$ such that $F\subset E$. If the infimum in \eqref{eq: LES} is $\lambda\in [1,\infty)$, then we say that $\X$ has $\lambda$-DES in direction $f$. 

Clearly finite-dimensional spaces have DES in every direction, since the unique minimum volume ellipsoid containing the unit ball is bounded. We do not know whether the property of having DES in all directions is inherited by 
the subspaces. The directionally Euclidean structure is preserved under the isometry group in the sense that if $\X$ has DES in direction 
$f\in \S_{\X^{\ast}}$ and $T\in \mathcal{G}_{\X}$, then $\X$ has DES in direction $T^{\ast}f$ also. On the other hand, DES is not preserved in isomorphisms, as it turns out. For instance, for each $x\in \S_{\X}$ one can present $\X$ isomorphically as $[x]\oplus_{2}\Y$ 
for suitable $\Y\subset\X$ and it turns out that $[x]\oplus_{2}\Y$ has DES in the direction of the axis $[x]$, roughly speaking.\ \\

Let us begin with an observation on the existence of continuous inner products with controlled spread in one direction. 
\ \\
\begin{theorem}\label{thm: 123}
Let $f\in\S_{\X^{\ast}}$ and consider the following conditions:\\
\begin{enumerate}
\item[(1)]{$\X$ has DES in direction $f$.}
\item[(2)]{$\sup_{F}\inf_{E}\eta(f,E)<\infty$, where $F\subset E$ are finite-dimensional}
\item[(3)]{There exists (a possibly degenerate) inner product $(\cdot|\cdot)$ on $\X$ such that 
$(x|y)\leq ||x||\cdot ||y||$ for $x,y\in \X$ and $(z|z)\geq \left(\frac{|f(z)|}{\sup_{F}\inf_{E}\eta(f,E)}\right)^{2}$ for each $z\in \X$.}
\end{enumerate}
Then $(1)\implies (2)\implies (3)$.
\end{theorem}

In the proof of this result we require the following combinatorial fact.
\begin{lemma}\label{lm: DS}
Let $(X,\leq)$ be a directed set, $(Y,\preceq)$ a poset and $f\colon X\rightarrow Y$ a map.
Suppose that there exists $y\in Y$ with the following property: For each $x\in \X$ there exists
$z\in X$ such that $x\leq z$ and $f(z)\preceq y$. Then there exists a subset $Z\subset X$ satisfying the following
conditions:
\begin{enumerate}
\item[(i)]{$(Z,\leq)$ is a directed set.}
\item[(ii)]{For each $x\in\X$ there is $z\in Z$ with $x\leq z$.}
\item[(iii)]{$f(z)\preceq y$ for each $z\in Z$.}
\end{enumerate}
\end{lemma}
\begin{proof}
Fix $y\in Y$ and $f$ as above. We denote the choice function $\alpha\colon \mathcal{P}(\X)\setminus\{\emptyset\}\rightarrow \X,\ \alpha(A)\in A,$ 
provided by Axiom of Choice. We will proceed by recursion of countable length as follows.

Initiate the recursion by putting
\[A_{0}=\{\alpha(\{z\in X|\ f(z)\preceq y,\ x\leq z\})|\ x\in \X\}.\]
Given $i<\omega$ and $A_{i}$ we define $A_{i+1}$ by
\[A_{i+1}=A_{i}\cup \{\alpha(\{z\in X|\ f(z)\preceq y,\ a,b\leq z\})|\ a,b\in A_{i}\}.\]
Then $Z=\bigcup_{i<\omega}A_{i}$ defines the claimed subset of $\X$.
\end{proof}  

\begin{proof}[Proof of Theorem \ref{thm: 123}]
Instead of proving $(1)\implies (2)$ we will prove a stronger statement, namely that $\sup_{F}\inf_{E}\eta(f,E)\leq \inf_{F}\sup_{E}\eta(f,E)$. Towards this, put $\epsilon>0$. 
Observe that if $F_{0}$ is a finite-dimensional subspace such that $\sup_{E_{0}}\eta(f,E)+\epsilon \leq \inf_{F}\sup_{E}\eta(f,E)$,
where $F_{0}\subset E_{0}$ are finite-dimensional, then 
\[\sup_{F}\inf_{E}\eta(f,E) +\epsilon\leq \inf_{F}\sup_{E}\eta(f,E)\] 
because in the left the infimum can be taken over $E$ such that $F,F_{0}\subset E$. Thus the claim holds as $\epsilon$ was arbitrary.

To check direction $(2)\implies (3)$, let $\mathcal{F}$ be the set of all finite-dimensional subspaces of $\X$ ordered by inclusion. Consider the map $\alpha\colon \mathcal{F}\to [0,\infty)$ given by $E\mapsto \eta(f,E)$. 

According to Lemma \ref{lm: DS} and the fact $\sup_{F}\inf_{E}\alpha(E)<\infty$ we obtain that there is for each $i\in\N$ a directed subset $\mathcal{F}_{i}\subset\mathcal{F}$ satisfying the statements $(i)$, $(ii)$ 
of the lemma and $\alpha(E)\leq \sup_{F}\inf_{E}\alpha(E)+\frac{1}{i}$ for $E\in \mathcal{F}_{i}$. Let $\mathcal{M}=\bigcup_{i}\mathcal{F}_{i}$. 

For each $E\in\mathcal{M}$ let $(\cdot|\cdot)_{E}$ be the inner product on $E$ corresponding to the 
minimum volume ellipsoid $\mathcal{E}_{E}$. For each $E\in \mathcal{M}$ let $P_{E}\colon \X\to E$ be a (bounded) linear projection.
Consider $\R^{\mathcal{M}}$ with the point-wise linear structure. Define a map $B\colon \X\times \X \to \R^{\mathcal{M}}$ by $B(x,y)(E)=(P_{E}x|P_{E}y)_{E}$. 
Clearly $B$ is a bilinear map. 

The family 
\[\{\{E\in\mathcal{M}|\ F\subset E,\ \alpha(E)\leq \inf_{E}\alpha(E)+i^{-1}\}\}_{(F,i)\in\mathcal{M}\times \N}\] 
is a filter base on $\mathcal{M}$.
Let $\mathcal{U}$ be an ultrafilter extending the above filter base. Put $(x|y)_{\X}=\lim_{E,\mathcal{U}}B(x,y)(E)$ for $x,y\in\X$. It is easy to see that $(\cdot|\cdot)_{\X}$ is a bilinear mapping. 

Let us check that $(x|y)_{\X}\leq ||x||\cdot ||y||$ for $x,y\in \X$. Indeed, let us first observe that the set $\{E\in \mathcal{M}:\ x,y\in E\}$ contains sets in the filter base and thus 
this set is in $\mathcal{U}$. For this reason $\lim_{E,\mathcal{U}}(P_{E}(x)-x)=0$ and $\lim_{E,\mathcal{U}}(P_{E}(y)-y)=0$ . 
Since $\B_{E}\subset \mathcal{E}_{E}$ for each $E$, we get that $(x|x)_{E}\leq ||x||^{2}$ for $x\in \X$ and $E\in\mathcal{M}$ such that $x\in E$. 
Thus 
\[\lim_{E,\mathcal{U}}(P_{E}x|P_{E}y)_{E}\leq \lim_{E,\mathcal{U}}\sqrt{(P_{E}x|P_{E}x)_{E}(P_{E}y|P_{E}y)_{E}}\leq \lim_{E,\mathcal{U}}\sqrt{||x||^{2}\cdot ||y||^{2}}=||x||\cdot ||y||.\]

Observe that according to the selection of the filter basis, we get that $\lim_{D,\mathcal{U}}(\alpha(D)-\sup_{F}\inf_{E}\eta(f,E))=0$. Let $\delta_{E}(x)=\sup\{a>0:\ ax\in \mathcal{E}_{E}\}$ for $x\in E,\ x\neq 0$, 
$E\in \mathcal{M}$. Since $|f(\delta_{E}(x)x)|\leq \alpha(E)$ for each $E$ such that $x\in E$, we obtain that
\[\frac{|f(x)|}{\alpha(E)}\leq \frac{1}{\delta_{E}(x)},\]
which yields 
\[\frac{|f(x)|}{\sup_{F}\inf_{E}\alpha(E)}=\lim_{E,\mathcal{U}}\frac{|f(x)|}{\alpha(E)}\leq \lim_{E,\mathcal{U}}\frac{1}{\delta_{E}(x)}=\sqrt{(x|x)_{\X}}.\]
\end{proof}

There is one considerable difference between the quantities '$\liminf$' and '$\limsup$' above, that is, $\sup_{F}\inf_{E}\eta(f,E)$ and $\inf_{F}\sup_{E}\eta(f,E)$. Namely, by using the former, weaker, control we may construct
inner products, which do not vanish in the given direction $f$. Later we wish to construct inner products with a control simultaneously in several directions, and then the latter, stronger, concept is required. Next, it turns out that DES in every direction 
is a strong enough property of Banach spaces to characterize Hilbert spaces up to isomorphism. 

\begin{theorem}\label{thm: main}
Let $\X$ be a Banach space. The following conditions are equivalent:
\begin{enumerate}
\item[(i)]{$\X$ is isomorphic to a Hilbert space.} 
\item[(ii)]{$\X$ is isomorphic to a space $\Y$ having has DES in every direction $f\in \S_{\Y^{\ast}}$.}
\end{enumerate}
Moreover, $\X$ is isometric to a Hilbert space if and only if it has $1$-DES in every direction $f\in \S_{\X^{\ast}}$.
\end{theorem}
\begin{proof}
First observe that a Hilbert space clearly has $1$-DES in every direction. This observation covers the implication $(i)\implies (ii)$ and the latter 'only if' part.

Suppose that $\X$ is isomorphic to $\Y$ having DES in all directions $f\in \S_{\X^{\ast}}$. 
Let us resume the notations of the proof of Theorem \ref{thm: 123}: we let $\mathcal{F}$ be the set of finite-dimensional subspaces of $\Y$ ordered by inclusion.
Let $\mathcal{U}$ here be an ultrafilter on $\mathcal{F}$ containing the filter basis $\{\{E\in \mathcal{F}:\ F\subset E\}\}_{F\in \mathcal{F}}$.
We still denote $(x|y)_{\Y}=\lim_{E,\mathcal{U}}(P_{E}x|P_{E}y)_{E}$. By the construction $(\cdot|\cdot)_{\Y}$ is a bilinear map with $(x|y)_{\Y}\leq ||x||\cdot ||y||$ for $x,y\in\X$.

Now let us study the 'ball' $B=\{y\in \Y:\ (y|y)_{\Y}\leq 1\}$. Pick $b\in B$ with $\lim_{E,\mathcal{U}}(P_{E}b|P_{E}b)_{E}< 1$. We obtain that 
\begin{equation}\label{eq: est}
|f(b)|\leq \lim_{E,\mathcal{U}}\sup_{c\in \mathcal{E}_{E}}|f(c)|\leq \inf_{F}\sup_{E}\eta(f,E).
\end{equation}
If $\Y$ has DES in all directions, then the right hand side of \eqref{eq: est} is finite. Thus we obtain that $\sup_{b\in B}|f(b)|<\infty$.
Now, since $f\in \S_{\Y^{\ast}}$ was arbitrary, the Uniform Boundedness Principle yields that $B$ is norm-bounded. Observe that $\B_{Y}\subset B$ by the selection of the ellipsoids $\mathcal{E}_{E}$.
This can be rephrased as follows: there exists $1\leq C<\infty$ such that 
\[||y||^{2}\leq (y|y)_{\Y}\leq C||y||^{2},\quad \mathrm{for}\ y\in \Y.\] 

To check the last claim, we will apply the above argument with $\X=\Y$, where this space has $1$-DES in all directions. 
It follows that $|f(b)|\leq 1$ in \eqref{eq: est}. Since this holds for all $b\in B$ and all $f\in \S_{\X^{\ast}}$, we get that $B\subset \B_{\X}$. Thus $B=\B_{\X}$ and hence $||x||^{2}=(x|x)_{\X}$ for $x\in\X$. 
\end{proof}
\begin{remark}\label{remark: norming}
We note that in Theorem \ref{thm: main} the condition $(ii)$ could be replaced by an equivalent condition $(ii^{\prime})$: 
$\X$ is isomorphic to a space $\Y$ and there exists a norming subspace $\Z\subset \Y^{\ast}$ such that $\Y$ has DES in every direction $f\in \S_{\Z}$.
\end{remark}

\subsection{On boundedness of ellipsoids containing the unit ball}

In a Banach space a continuous non-degenerate inner product corresponds to a convex body that might be unbounded but does not contain a $1$-dimensional linear subspace. There are plenty of such convex bodies, or 
inner products, according to the following fact: 

\begin{proposition}
Let $\X$ be a Banach space with a $\omega^{\ast}$-separable dual space. Then there exists an inner product $(\cdot|\cdot)$ on $\X$ such that 
$(x|x)\leq ||x||^{2}$ and $(x|x)>0$ for $x\in\X,\ x\neq 0$. This in turn implies that there exists a continuous linear injection from $\X$ into $c_{0}(\Gamma)$ for some set $\Gamma$. 
\end{proposition}
\begin{proof}
Let $(x_{n}^{\ast})\subset \X^{\ast}$ be a sequence, which is $\omega^{\ast}$-dense in $\X^{\ast}$ and does not contain $0$. It is easy to see that that the mapping 
$g\colon \X\to \ell^{\infty}$ given by $x\mapsto (\frac{x_{n}^{\ast}(x)}{||x_{n}^{\ast}||})$ is linear, contractive and injective. The latter property follows from the fact that 
$(x_{n}^{\ast})$ separates each $x\in \X,\ x\neq 0$, as it is $\omega^{\ast}$-dense. Next, observe that $f\colon \ell^{\infty}\to \ell^{2}$ given by 
$(z_{n})\mapsto (2^{-n}z_{n})$ is linear, contractive and injective. Now, the required inner product on $\X$ is induced by one in $\ell^{2}$ via the composite mapping $f\circ g$.

In the latter claim we first form the completion $\H$ of $(\X,||\cdot||_{2})$ and consider its orthonormal basis $\{e_{\gamma}\}_{\gamma\in\Gamma}$. The required map 
$\X\to c_{0}(\Gamma)$ is then given by $x\mapsto \{(x|e_{\gamma})\}_{\gamma\in\Gamma}$.
\end{proof}
Spaces that admit a continuous linear injection into $c_{0}(\Gamma)$ have been studied quite a bit (see e.g. \cite{God}). We do not know exactly what kind of Banach spaces admit a continuous (non-degenerate) inner product.
Note that a continuous inner product will induce a continuous norm on $\X$, but according to the Open Mapping Principle this norm will be \emph{complete} if and only if the original norm and the induced norm are equivalent.

Even though there are often continuous inner products on Banach spaces, the way the inner products are produced here by applying DES results in constructions that are sensitive to small changes of the norm, as it turns out. 
We also note that continuous inner products are typically not invariant under isometries of the space. However, directionally Euclidean structure allows us to construct continuous inner products that will be invariant under 
suitable isometries.

The following result suggests that the disposition of the minimum volume ellipsoids becomes in a sense chaotic as the subspaces vary. It follows in particular that the property of having DES in all
directions is not preserved in isomorphisms, even for small Banach-Mazur distances.

\begin{theorem}
Let $\X$ be a Banach space and $F,E\subset \X$ be subspaces such that $\X=F\oplus_{p} E$ isometrically for some $1\leq p\leq \infty$. If $1\leq p\leq 2$ and $\dim(F)<\infty$, 
then $\X$ has DES in every direction $(f,0)\in \S_{F^{\ast}\oplus_{p^{\ast}}E^{\ast}}$. On the other hand, if $2<p\leq \infty$ and $E=\ell^{2}$, then $\X$ does not have DES in any direction $(f,0)\in \S_{F^{\ast}\oplus_{p^{\ast}}\ell^{2}}$.
\end{theorem}

\begin{proof}
For both the cases $p\leq 2$ and $p>2$ we are interested in finite-dimensional subspaces of the type $F\oplus_{p} E_{n}$, where $E_{n}\subset E$ is an $n$-dimensional subspace. 
This is so because each finite-dimensional subspace $\Y\subset \X$ is contained in a finite-dimensional subspace of the above type. 

Denote by $\mathcal{E}\subset F\oplus_{p} E_{n}$ the unique minimal volume ellipsoid containing $\B_{F\oplus_{p} E_{n}}$.
Let $(\cdot|\cdot)_{\mathcal{E}}$ be the corresponding inner product. According to the results by Auerbach this ellipsoid is invariant under linear isometries of 
$F\oplus_{p} E_{n}$ onto itself. In particular, given a linear projection $P_{F}\colon F\oplus_{p} E_{n}\to F$ it holds that $\mathcal{E}$ is invariant under the isometric reflection mapping
$\I-2P_{F}$. Denote $Q_{F}=\I-P_{F}$.

\textit{Claim 1.} For each $x\in F\oplus_{p} E_{n}$ it holds that 
\[(x|x)_{\mathcal{E}}=(P_{F}x|P_{F}x)_{\mathcal{E}}+(Q_{F}x|Q_{F}x)_{\mathcal{E}}.\]
Indeed, by using the invariance we obtain
\[(x|x)_{\mathcal{E}}=((\I-2P_{F})x|(\I-2P_{F})x)=(x|x)_{\mathcal{E}}-2(x|P_{F}x)_{\mathcal{E}}-2(P_{F}x|x)_{\mathcal{E}}+4(P_{F}x|P_{F}x)_{\mathcal{E}},\]
which gives that $(x|P_{F}x)_{\mathcal{E}}=(P_{F}x|P_{F}x)_{\mathcal{E}}$, and this yields the claim.

Next, recall that each ellipsoid is given by the formula $\frac{x_{1}^{2}}{C_{1}}+\frac{x_{2}^{2}}{C_{2}}+\ldots +\frac{x_{m}^{2}}{C_{m}}\leq 1$ where 
one fixes a suitable coordinate system. If one normalizes the measure by fixing the volume of $\B_{F\oplus_{p} E_{n}} $, then the volume of the ellipsoid does not depend on 
the selection of the coordinate system. With the above notation the volume of the ellipsoid $\mathcal{E}$ is 
\[\mathrm{Vol}(\mathcal{E})=\beta C_{1}C_{2}\dots C_{m},\]
where $\beta$ is a constant depending on the dimension and the coordinate system. 

Let us identify $F\oplus_{p} E_{n}$ with $\R^{m}$. Thus we will regard the volume as the standard $m$-dimensional Lebesgue measure. According to Claim 1. we may assume here that the coordinate system 
is given in such a way that the first $\dim(F)$-many coordinates of $\R^{m}$ support the ellipsoid $P_{F}(\mathcal{E})$ and the last $n$ coordinates 
support the ellipsoid $Q_{F}(\mathcal{E})$. We may assume without loss of generality, by choosing the coordinate system suitably, that $\B_{\ell^{2}(k+n)}\subset \B_{F\oplus_{p}E}$. 
Clearly, the fact that $\B_{\ell^{2}(k)}\subset P(\mathcal{E})$ and $\B_{\ell^{2}(n)}\subset Q(\mathcal{E})$ implies that $C_{1},C_{2},\ldots, C_{k+n}\geq 1$.

Fix the minimal volume ellipsoids $\mathcal{E}_{1}\subset F$, $\mathcal{E}_{2}\subset E_{n}$ containing $\B_{F}$ and $\B_{E_{n}}$, respectively.
Then the corresponding constants $C^{(F)}_{1}C^{(F)}_{2}\dots C^{(F)}_{k}\leq C_{1}C_{2}\dots C_{k}$ and $C^{(E)}_{k+1}C^{(E)}_{k+2}\dots C^{(E)}_{m}\leq C_{k+1}C_{k+2}\dots C_{m}$.

Let us verify the statement involving $p\leq 2$. We obtain that 
\[\B_{F\oplus_{p} E_{n}}\subset\{x\in F\oplus_{p} E_{n}:\ (P_{F}x|P_{F}x)_{\mathcal{E}_{1}}+(Q_{F}x|Q_{F}x)_{\mathcal{E}_{2}}\leq 1\}.\]
We conclude that 
\[(x|x)_{\mathcal{E}}= (P_{F}x|P_{F}x)_{\mathcal{E}_{1}}+(Q_{F}x|Q_{F}x)_{\mathcal{E}_{2}}\quad \mathrm{for}\ x\in F\oplus_{p} E_{n}.\]
This means that $P_{F}(\mathcal{E})$ is a norm bounded set, which does not depend on $n=\dim(E_{n})$. Consequently we have the first part of the statement.

Let us check the latter statement, where $p>2$ and $E$ is a Hilbert space. The invariance of $\mathcal{E}$ under isometries yields that $C_{k+1}=C_{k+2}=\ldots =C_{m}$. 
Indeed, here we consider isometries of the form $\I\oplus T$, where $\I\colon F\to F$ is the identity map, $T$ is a linear isometry of $E_{n}=\ell^{2}(n)$ onto itself, and 
we apply the fact that the isometry group of a Hilbert space acts transitively on the unit sphere. 

\textit{Claim 2.} Given $b>1$ and $2<p<\infty$ we denote
\[\alpha_{p}(b)=\inf\{a>1:\ \left(\frac{x^{2}}{a}+\frac{y^{2}}{b}\right)^{\frac{1}{2}}\leq (x^{p}+y^{p})^{\frac{1}{p}}\ \mathrm{for}\ x,y>0\}.\]
Then by analyzing with the test points $(2^{-\frac{1}{p}},2^{-\frac{1}{p}})$ and $(1,1)$ we obtain the following estimate:
\begin{equation}\label{eq: apbb}
\frac{b}{2^{\frac{2}{p}}b-1}\leq \alpha_{p}(b)\leq \frac{b}{b-1}.
\end{equation}
Moreover, it is fairly easy to see that 
\begin{equation}\label{eq: ainfty}
\alpha_{p}(b)\to \infty\ \mathrm{as}\ b\to 1^{+}\quad \mathrm{for}\ 2<p\leq \infty .
\end{equation}

Note that according to Claim 2. the constants $C_{1},C_{2},\dots, C_{k}$
satisfy the inequality 
\begin{equation}\label{eq: apb}
a\frac{C_{k+1}}{2^{\frac{2}{p}}C_{k+1}-1}\leq C_{i}\leq b\frac{C_{k+1}}{C_{k+1}-1},\quad 1\leq i\leq k
\end{equation}
for suitable constants $a,b>0$ depending only on the disposition of $\B_{F}$ in $\R^{k}$, and not on $p$ or the actual value of $k$ or $n$.

Since $C_{k+1}=C_{k+2}=\cdots =C_{m}$, it follows by using \eqref{eq: apb}  that the expression of the minimal volume
\[\mathrm{Vol}(\mathcal{E})=\beta^{(n)}C_{1}^{(n)}C_{2}^{(n)}\dots C_{k}^{(n)}(C_{k+1}^{(n)})^{n},\]
must satisfy that $C_{k+1}^{(n)}\longrightarrow 1^{+}$ as $n\to \infty$. Thus we obtain by \eqref{eq: ainfty} that $C_{i}^{(n)}\to \infty$ as $n\to \infty$ for $1\leq i\leq k$.
This yields the latter claim.
\end{proof}

We do not know whether the property of having DES in all (or some) directions passes on to subspaces, ultrapowers,  ultra-roots, or to almost isometric copies. 
Consider Banach spaces $\X$ with the property that the minimum volume ellipsoids $\mathcal{E}$ of finite-dimensional subspaces of $\X$ are uniformly bounded. 
We note that this property passes on to the abovementioned structures related to $\X$. The proof is omitted but we will list some observations, which lead to the claims involving ultrapowers.
\medskip \\
\noindent \textit{Observation 1: The determination of the volume of the minimum volume ellipsoid is continuous with respect to the norm.} This can be formulated more precisely as follows: suppose that $E$ and $F$ are $n$-dimensional
spaces with Lebesgue measures $\mu$ and $\nu$, respectively, normalized such that $\mu(\B_{E})=\nu(\B_{F})=1$. If $f\colon E\to F$ is an isomorphism and $||x||\leq ||f(x)||\leq C||x||$ for $x\in E$, then 
$\mu(A)\leq \nu(f(A))\leq C^{n}\mu(A)$ for each Lebesgue measurable $A\subset E$.
\medskip \\
\noindent \textit{Observation 2: Each finite-dimensional subspace of $\X^{\mathcal{U}}$ is an ultralimit of suitable finite-dimensional subspaces of $\X$ in the sense of the Banach-Mazur distance.} Let $x_{1},x_{2},\ldots,x_{n}\in \X^{\mathcal{U}}$
be non-zero, linearly independent vectors. Then for each $1\leq i\leq n$ there is a sequence $(z_{k}^{(i)})\subset \X$ such that $\lim_{k,\mathcal{U}}||z_{k}^{(i)}-x_{i}||=0$. Let $E_{k}=[z_{k}^{(1)},z_{k}^{(2)},\ldots,z_{k}^{(n)}]$ for 
$k\in \N$. It is not difficult to see that there exists $K\in \mathcal{U}$ such that $T_{k}\colon E_{k}\to E,\ T_{k}\left(\sum_{i} a_{i}z_{k}^{(i)}\right)=\sum_{i} a_{i}x_{i}$ defines a linear isomorphism for $k\in K$. Moreover, 
$\lim_{k,\mathcal{U}}||T_{k}||=1$ and $\lim_{k,\mathcal{U}}||T_{k}^{-1}||=1$.
\medskip \\
For each $k$ we denote by $\mu_{k}$ the Lebesgue measure on $E_{k}$ normalized such that $\mu_{k}(E_{k})=1$. We will denote the minimum volume ellipsoid of $E_{k}$ containing $\B_{E_{k}}$ by $\mathcal{E}_{k}$. The minimum volume ellipsoid
of $E$ is denoted by $\mathcal{E}$ and the corresponding normalized measure is denoted by $\mu$.
\medskip \\
\noindent \textit{Observation 3.} By using the previous observations we obtain that $\lim_{k,\mathcal{U}}\mu_{k}(\mathcal{E}_{k})=\mu(\mathcal{E})$.
\medskip \\
\noindent \textit{Observation 4.} By using the previous observations, compactness and the uniqueness of the minimum volume ellipsoid in $E$, 
we have that $\lim_{k,\mathcal{U}}(T_{k}^{-1}x|T_{k}^{-1}y)_{\mathcal{E}_{k}} = (x|y)_{\mathcal{E}}$ for $x,y\in E$.

The details are omitted.

\section{An application involving Mazur's problem}

In the presence of a high degree of symmetry the nice local properties of a Banach space tend to self-improve, see e.g. \cite{Fi}, also \cite{BR}. For example, a convex-transitive Banach space with the RNP is already
uniformly convex, uniformly smooth and almost transitive.

In \cite{Ca} it was asked whether a Banach space almost transitive with respect to isometric finite-dimensional perturbations of the identity is isometric to a Hilbert space. 
The answer is affirmative in the case that the space has DES at some direction.

\begin{theorem}
Let $\X$ be a Banach space, which has DES at some direction $f\in \S_{\X}$ and assume that $\X$ is convex-transitive with respect to $\mathcal{G}_{F}$. 
Then $\X$ is isometrically a Hilbert space.
\end{theorem} 
\begin{proof}
We will construct an inner product $(\cdot|\cdot)_{\X}$ on $\X$ such that $|x|=\sqrt{(x|x)_{\X}}$ defines 
a norm, which is continuous with respect to $||\cdot||$ and such that $|x|=|Tx|$ for $x\in \X$ and $T\in \mathcal{G}_{F}$.
By using that $\X$ is convex-transitive with respect to $\mathcal{G}_{F}$ it follows is that there exists
a constant $c>0$ such that $||\cdot||=c|\cdot|$ (see \cite{Co}), which yields the claim. 

The argument here closely resembles that of \cite{Ca}. Let $\Gamma$ be the set of all finite subsets $\gamma$ of $\mathcal{G}_{F}$ such that $T\in \gamma\implies T^{-1}\in \gamma$.
Note that $\Gamma$ can be viewed as a lattice when ordered by inclusion $\subseteq$.

By using a simple argument concerning the Hamel basis of $\X$ we obtain that for each finite-dimensional subspace $A\subset \X$ and each $\gamma\in \Gamma$ there is a finite-dimensional subspace $F\supset A$ and 
a finite-codimensional subspace $E$ such that $\X=F\oplus E$, $\span\bigcup_{T\in \gamma}(\I-T)(\X)\subset F$, 
$E\subset \bigcap_{T\in \gamma}\ker(\I-T)$ and $T(F)=F$ for $T\in \gamma\in \Gamma$.
Note that $F\neq \{0\}$ for any $\gamma\neq \{\I\}$. Denote by $\mathcal{F}_{A,\gamma}$ the collection of all pairs $(F,\gamma)$, which satisfy the above conditions.

For each $(F,\gamma)\in \mathcal{F}_{A,\gamma}$ let $\mathcal{E}_{F}$ be the unique minimum-volume ellipsoid in $F$, which contains
$\B_{\X}\cap F$. We denote by $(\cdot|\cdot)_{F}\colon F\times F\to \R$ the inner product induced by the 
ellipsoid $\mathcal{E}_{F}$. Observe that the uniqueness of the minimum value ellipsoid implies invariance under isometries and thus
$(Tx|Ty)_{F}\leq (Tx|Tx)_{F}(Ty|Ty)_{F}=(x|x)_{F}(y|y)_{F}\leq ||x||\ ||y||$ for $x,y\in F,\ T\in\gamma$.

For technical reasons, for each finite-dimensional $F\subset \X$ we denote a linear projection $\X\to F$ by $P_{F}$ without specifying exactly which projection is in question. 
Let $\mathcal{M}=\bigcup_{A,\gamma}\mathcal{F}_{A,\gamma}$, where the union is taken over finite-dimensional subspaces $A\subset \X$ and $\gamma\in \Gamma$. 
We may define a partial order on $\mathcal{M}$ by declaring $(F,\gamma)\leq (F^{\prime},\gamma^{\prime})$ if $F\subset F^{\prime}$ and $\gamma\subset \gamma^{\prime}$.

Define $[\cdot|\cdot]\colon \X\to \R^{\mathcal{M}}$ by letting $[x|y]$, evaluated at $(F,\gamma)$, equal to $(P_{F}x|P_{F}y)_{F}$. 
Observe that the family 
\[\{\{(F,\gamma)\in \mathcal{M}:\ \delta\subset\gamma,\ A\subset F\}\}_{(A,\delta)\in \mathcal{M}},\]
where $A$ ranges in finite-dimensional spaces, defines a filter basis on $\mathcal{M}$. 

Let $\mathcal{U}$ be a non-principal ultrafilter on $\mathcal{M}$ extending this filter basis.
Define $(\cdot|\cdot)\colon \X\times \X\rightarrow \R$ by $(x|y)=\lim_{\mathcal{U}}[x|y]$.
It is easy to check to that $(\cdot|\cdot)$ is well-defined, bilinear and $(x|y)\leq ||x||\ ||y||$, $(x|x)\geq 0$, $(Tx|Ty)=(x|y)$ for $x,y\in E,\ T\in \mathcal{G}_{F}$.  
Indeed, pick $(F,\gamma)\in \mathcal{M}$ such that $x,y\in F,\ T\in \gamma$. Firstly, $[\cdot|\cdot]$ evaluated at $(F^{\prime},\gamma^{\prime})\geq (F,\gamma)$ satisfies the conditions described above. 
Secondly, note that the set of pairs $(F^{\prime},\gamma^{\prime})\geq (F,\gamma)$ belongs to $\mathcal{U}$. This means that the ultralimit $(\cdot|\cdot)=\lim_{\mathcal{U}}[\cdot|\cdot]$ satisfies the abovementioned conditions.

Finally, we will check that $(x|x)>0$ for $x\in \X,\ x\neq 0$. Suppose that $\X$ has $\lambda$-DES at direction $f\in \S_{\X^{\ast}}$. Then $\inf_{F}\sup_{F^{\prime}}\eta(f,F^{\prime})<\infty$, which means that 
we may select pick $F$ such that $\sup_{F^{\prime}\supset F}\eta(f,F^{\prime})=\alpha<\infty$. Pick $x\in \X$ such that $f(x)\geq \alpha$. Then, $[x|x]$ evaluated at any pair $(F^{\prime},\gamma^{\prime})\geq (F,\gamma)$ 
has value at least $1$. Similarly as above, since the set of pairs $(F^{\prime},\gamma^{\prime})\geq (F,\gamma)$ belongs to the filter, we obtain that the ultralimit $(x|x)\geq 1$.
This means that $\Y=\{x\in \X:\ (x|x)=0\}$ is not the whole space $\X$. Observe that $\Y$ is invariant under $\mathcal{G}_{F}$. It is easy to see that the convex-transitivity with respect to $\mathcal{G}_{F}$ yields
that only a trivial subspace, i.e. $\{0\}$ or $\X$, can be invariant under $\mathcal{G}_{F}$. Since $\Y\neq \X$, we conclude that $\Y=\{0\}$. 
\end{proof}

\begin{theorem}
Let $\X$ be a Banach space, which is convex-transitive with respect to $\mathcal{G}_{\X}$ and has DES in some direction $f\in \S_{\X^{\ast}}$. Then $\X$ is isomorphic to a Hilbert space.
Moreover, if $\X$ has $1$-DES in direction $f$, then $\X$ is isometric to a Hilbert space.
\end{theorem}
\begin{proof}
Suppose that $\X$ has $\lambda$-DES in direction $f\in \S_{\X^{\ast}}$.
It is known (see \cite{BR}) that $\X$ is convex-transitive if and only if $\overline{\conv}^{\omega^{\ast}}(\{T^{\ast}g:\ T\in \mathcal{G}_{\X}\})=\B_{\X^{\ast}}$ for $g\in\S_{\X^{\ast}}$. 
This means that $\{T^{\ast}f:\ T\in \mathcal{G}_{\X}\}$ is a $1$-norming set. Then one can construct, similarly as in the proof of Theorem \ref{thm: main}, an inner product $(\cdot|\cdot)_{\X}$ on $\X$, 
which satisfies that $(x|x)_{\X}\leq ||x||^{2}$ for $x\in \X$. It follows from the assumptions by inspecting the construction of $(\cdot|\cdot)_{\X}$ that $\{x\in\X:\ (x|x)_{\X}\leq 1\}\subset \lambda\B_{\X}$. 
Thus we have the claim.
\end{proof}

Some authors have asked whether an almost transitive Banach space, that is isomorphic to a Hilbert space, is in fact isometric to one (see e.g. \cite{Ca0}). This question appears not to have been settled at the moment.

\section{Final remark: near-convexity of the duality mapping}

Recall that the duality mapping $J\colon \X\rightarrow 2^{\X^{\ast}}$ is a multivalued mapping defined by 
$J(x)=\{x^{\ast}\in \X^{\ast}:\ ||x||^{2}=||x^{\ast}||^{2}=x^{\ast}(x)\}$ for $x\in \X$.
If $\X$ is a Gateaux-smooth space, then $J$ becomes a point-to-point mapping. Recall that for Hilbert spaces the duality map is an isometric isomorphism.
We always have that $\overline{J(\B_{\X})}=\B_{\X^{\ast}}$ according to the Bishop-Phelps theorem and in the reflexive case $J(\B_{\X})=\B_{\X^{\ast}}$ according to James's characterization of reflexivity.  However, 
it easily happens that the image of a convex set under the mapping $J$ is not convex. Next, we will study spaces whose duality mapping does not distort convex sets very far from being convex. It turns out that such spaces are isomorphically Hilbertian.

\begin{theorem}
Let $\X$ be a smooth Banach space and let $J\colon \X\to \X^{\ast}$ be the duality mapping. Suppose that there exists a constant
$0\leq C<1$ such that 
\[\left\|J\left(\sum x_{n}\right)-\sum J(x_{n})\right\|\leq C\left\|\sum x_{n}\right\|\quad \mathrm{for}\ x_{1},\ldots,x_{n}\in \X.\]
Then $\X$ is isomorphic to a Hilbert space. Moreover, if $J$ is a convex map, then $\X$ is isometric to a Hilbert space.
\end{theorem}
\begin{proof}
By using an ultrafilter construction similar to that in the proof of Theorem \ref{thm: main} it suffices to check that in any finite-dimensional subspace $F\subset \X$ 
there exists an inner product $(\cdot|\cdot)\colon F^{2}\to \R$ such that 
\begin{equation}\label{eq: (1-C)}
(1-C)||x||^{2}\leq (x|x)\leq (1+C)||x||\quad \mathrm{for}\ x\in F. 
\end{equation}

Let $F\subset \X$ be a finite-dimensional subspace. Recall that there exists an Auerbach basis on $F$, that is, a 
biorthogonal system $\{(e_{i},e_{i}^{\ast})\}_{i=1}^{n}\in (\S_{\X}\times \S_{\X^{\ast}})^{n}$. 
Define a mapping $g\colon F^{2}\to \R$ by $g(x,y)=\sum a_{i} e_{i}^{\ast}(x)$, where 
$\sum a_{i}e_{i}$ is the unique expression of $y$. Note that $g$ is bilinear. Define $B\colon F^{2}\to \R$ by $B(x,y)=\frac{g(x,y)+g(y,x)}{2}$
for $x,y\in F$. Now $B$ is clearly a symmetric bilinear form.

Fix $x=\sum a_{i}e_{i}\in F$. Observe that 
\begin{equation*}
\begin{array}{lll}
|\ ||x||^{2}-g(x,x)\ |&=&|J(x)(x)-g(x,x)|=\left\|\left(J\left(\sum a_{i}e_{i}\right)-\sum a_{i}J(e_{i})\right)(x)\right\|\\
&\leq & C ||x|| \cdot ||x||.
\end{array}
\end{equation*}
Thus $|\ ||x||^{2}-B(x,x)\ |\leq C||x||^{2}$, so that $B$ is a non-degenerate inner product on $F$ and \eqref{eq: (1-C)} holds. 
\end{proof}

\end{document}